\documentclass[12 pt]{amsart}
\usepackage{amsmath,amssymb,amsthm,amsfonts,setspace,hyperref,dsfont,yhmath,euscript}

\usepackage{color}

\newcommand{\ZZ}{\mathbb{Z}}

\newtheorem{theorem}{Theorem}[section]
\newtheorem{corollary}[theorem]{Corollary}

\newtheorem{lemma}[theorem]{Lemma}
\newtheorem{proposition}[theorem]{Proposition}

\newcommand{\comment}[1]{}
\theoremstyle{definition}

\numberwithin{equation}{section}

\numberwithin{equation}{section}

      \makeatletter
      \def\@setcopyright{}
      \def\serieslogo@{}
      \makeatother

\def\c{\mathcal{C}}
\def\R{\mathbb R}
\def\Q{\mathbb Q}

\def\Z{\mathbb Z}
\def\N{\mathbb N}
\def\T{\mathbb T}

\def\Id{\text{Id}}
\def\e{\epsilon}

\def\b{\beta}
\def\la{\lambda}

\def\A{\mathcal{A}}
\def\B{\mathcal{B}}

\def\E{{E}}

\def\M{\mathcal{M}}

\makeatletter
\renewcommand*\env@matrix[1][*\c@MaxMatrixCols c]{%
  \hskip -\arraycolsep
  \let\@ifnextchar\new@ifnextchar
  \array{#1}}
\makeatother
\begin{document}
\title[Local rigidity]
{Local rigidity for hyperbolic toral automorphisms}

\author[Boris Kalinin$^1$ \and Victoria Sadovskaya$^2$ \and Zhenqi Jenny Wang$^3$ ]{Boris Kalinin$^1$ \and Victoria Sadovskaya$^2$ \and Zhenqi Jenny Wang$^3$ }


\address{Department of Mathematics, The Pennsylvania State University, University Park, PA 16802, USA.}
\email{kalinin@psu.edu, sadovskaya@psu.edu}

\address{Department of Mathematics\\
        Michigan State University\\
        East Lansing, MI 48824,USA}
\email{wangzq@math.msu.edu}


\thanks{$^1$  Supported in part by Simons Foundation grant 855238}
\thanks{$^2$ Supported in part by NSF grant DMS-1764216}
\thanks{$^3$ Supported in part by NSF grant DMS-1845416}

\begin{abstract}
We consider a hyperbolic toral automorphism $L$ and its  $C^1$-small perturbation $f$.
It is well-known that $f$ is Anosov and topologically conjugate to $L$, but a conjugacy $H$ is only H\"older continuous in general. We discuss conditions for  smoothness of $H$, such as conjugacy of the periodic data of $f$ and $L$, coincidence of their Lyapunov exponents, and weaker regularity of $H$, and\,  we summarize questions, results, and techniques  in this area. Then  we introduce our new results:  if $H$ is weakly differentiable then it is $C^{1+\text{H\"older}}$ and, if $L$ is also weakly irreducible, then $H$ is $C^\infty$.

\end{abstract}

\maketitle

\section{Introduction}

Hyperbolic automorphisms of tori are the prime examples of hyperbolic dynamical systems.
The action of a matrix $L \in SL(d,\Z)$ on $\R^d$ induces an automorphism of the torus
$\T^d=\R^d/\Z^d$, which we denote by the same letter. An automorphism $L$ is called {\em hyperbolic}
or {\em Anosov} if the matrix has no eigenvalues on the unit circle. In general, 
a diffeomorphism $f$ of a compact manifold $\M$ is called Anosov if there exist 
a continuous $Df$-invariant splitting $T\M=E^s\oplus E^u$ and constants $K>0$ and 
$\theta<1$ such that for all $n\in \N$,
$$
\| Df^n(v) \|  \leq K\theta^n \| v \|
     \;\text{ for all }v \in E^s , \;\text{ and }\;\, \| Df^{-n}(v) \| \leq K\theta^n  \| v \|
     \;\text{ for all }v \in E^u. 
$$  
\noindent
The sub-bundles $E^s$ and $E^u$ are called stable and unstable. They are tangent to the stable and unstable foliations $W^s$ and $W^u$. A diffeomorphism is {\it transitive} if there is a point in $\M$ with dense orbit. All known examples of Anosov diffeomorphisms have this property.
\vskip.1cm
One of the key properties of hyperbolic systems is {\it structural stability}: any diffeomorphism $f$ of $\M$ sufficiently $C^1$ close to an Anosov diffeomorphism 
$g$ is also Anosov and is topologically conjugate to $g$ \cite{A}. The latter means that there exists a ho\-meo\-morphism $H$ of $\M$, called a {\em conjugacy}, such that $g\circ H= H \circ f$. Moreover, $H$ is unique 
in a $C^0$ neighborhood of the identity. A conjugacy $H$ is always bi-H\"older, but it is usually not even $C^1$, as there are various obstructions to smoothness. 
The problem of establishing smoothness of the conjugacy $H$ from some weaker assumptions
has been extensively studied. It is often referred to as {\it local rigidity}, in the sense that
weak equivalence of $f$ and $g$ implies strong equivalence.

We will concentrate on a conjugacy $H$ between a hyperbolic toral automorphism $L$ and its small perturbation $f$. We have
\begin{equation} \label{IntrConj}
L\circ H= H \circ f, \quad\text{ that is, }\quad f= H^{-1} \circ L \circ{H}.
\end{equation}
Any two such conjugacies differ by an affine automorphism of $\T^d$ commuting with $L$ \cite{W},
and in particular have the same regularity. 

\section{Smoothness of the conjugacy}

\subsection{Necessary conditions for  smoothness of $H$.} 
Hyperbolic systems have abundance of periodic points, and a lot of information about the system is captured by the periodic data. If $f^n(p)=p$ for some $n\in \N$ then $L^n(H(p))=H(p)$. If $H$ is a $C^1$ diffeomorphism, then
differentiating the iterated conjugacy equation $f^n= H^{-1} \circ L^n \circ{H}$ at $p$ yields
$$
D_p(f^n)\,= \,(D_pH)^{-1}\circ L^n\circ D_pH\,=\,  C_p^{-1} \circ L^n \circ C_p, \quad\text{where}\quad C_p=D_pH,
$$ 
and so $f$ and $L$ have {\it conjugate periodic data}. This condition can be easily destroyed by a small perturbation near a fixed or periodic point, and hence $H$ is not $C^1$ in general.  
\vskip.1cm

 Let $\rho_i$ be the absolute values  of the eigenvalues of $L$, and let $E^i$ be the sum  of generalized eigenspaces corresponding to the eigenvalues with modulus $\rho_i$. We recall that the Lyapunov exponents of $L$ are $\log \rho_i$, more specifically, for every vector $v\ne 0$ in $E^i$, 
$$
 \lim_{n\to \pm \infty} \frac1n \log\| L^n (v)\|=\log \rho_i.
$$
If $H$ is a $C^1$ diffeomorphism, then the Lyapunov exponents of $f$ are the same as for $L$: for each $x\in \T^d$ and each $v$ with $ D_xH (v) \in E^i$, the Lyapunov exponent of $v$ is $\rho_i$.
We note that conjugacy of periodic data implies that the Lyapunov exponents of $f$ at each periodic point $p$ are the same as for $L$.  By periodic approximation \cite{WS,K11}, it follows that the Lyapunov exponents of $f$ for each invariant measure are also the same as for $L$.
\vskip.1cm

It is natural to ask whether necessary conditions, such as conjugacy of periodic data or 
equality of Lyapunov exponents for some or all measures, are also sufficient for $C^1$ regularity of $H$. 
In addition, one can consider sufficiency of weaker regularity properties of $H$ such  Lipschitz continuity, weak differentiability, and absolute continuity (in dimension two). Further, one can consider sufficiency 
of any of these conditions for  $C^\infty$ regularity of $H$.
\vskip.2cm 


\subsection{Two-dimensional case} 
Definitive results for Anosov diffeomorphisms of $\T^2$ were obtained by de la Llave and Moriy\'on in \cite{L0,LM, L1}. In this case, conjugacy of the periodic data, or equality of Lyapounov exponents at periodic points, or absolute continuity imply smoothness of $H$. More precisely, the following holds.

\begin{theorem} \cite{L1}\label{mainL1}
Let $f$ and $g$ be $C^k$, $k = 2, 3, \dots \infty, \omega$, Anosov diffeomorphism of $\T^2$ topologically conjugate by $H$. Then $H$ is $C^{k-\e}$ provided that either
\vskip.1cm
(a) the Lyapounov exponents of $f$ and $g$ at corresponding periodic orbits are the same 
\vskip.1cm
(b) or both $H$ and $H^{-1}$ are absolutely continuous.

\end{theorem}

\vskip.3cm

In most local rigidity results, the smoothness of $H$  is  first established along natural invariant foliations, 
such as the stable and unstable, or their finer sub-foliations. Global smoothness of $H$ is then obtained 
using regularity results such as Journ\'e lemma \cite{J}. 

The main approach in Theorem \ref{mainL1} is to show that $H$ maps the invariant volume, or more generally the SRB measure $\mu$, for $f$ to that for $g$.  For this either the periodic data or absolute continuity of $H$ can be used. Then $H$ maps the absolutely continuous conditional measures of $\mu$ along the leaves of $W^{u,f}$ to those for $W^{u,g}$. The densities of these measures on the leaves are positive H\"older continuous functions smooth along the leaves. Since the leaves are one-dimensional, any $H$ that matches the conditional measures is similarly smooth along the leaves. Smoothness along the stable foliations can be
obtained similarly. The study of regularity of these densities and related objects
was a major part of the proof of Theorem \ref{mainL1}.

Next, local rigidity results were obtained for systems conformal on full stable and unstable sub-bundles and 
for systems close to irreducible automorphisms with simple Lyapunov spectrum.

\subsection{Higher dimensional conformal case}

If $g$ is conformal on full stable and unstable sub-bundles,  smoothness of the conjugacy 
is given by the next theorem, provided that the perturbed system is also conformal.

\begin{theorem} \cite{KS03,L3} \label{conf}
 Let $f,g:\M\to \M$ be transitive $C^{k}$, $k = 2, 3, \dots \infty$, Anosov diffeomorphisms that are conformal (or more generally  uniformly quasiconformal) on their stable and unstable sub-bundles,
which have dimension at least two.  Than any topological conjugacy between $f$ and $g$ is $C^{k-\e}$.
\end{theorem}

For $k=\infty$ this was proved in \cite{KS03} under an additional assumption that either $\dim E^s$ and $\dim E^u $ are at least three or that $\M$ is a torus (or more generally, an infranilmanifold). It is a corollary of a global rigidity result: any such $f$ is $C^\infty$ conjugate to a finite factor of an automorphism of a torus \cite{KS03}. The result in \cite{L3} is for any $k$,  but it makes an additional assumption of conformality 
with respect to a $C^1$ conformal structure. In \cite{KS03,S05}, existence of a continuous invariant 
conformal structure is obtained from uniform quasiconformality, and its $C^\infty$ smoothness is then 
proved for $C^\infty$ $f$  using conformality and smoothness of stable and unstable holonomies. To verify the extra assumption \cite{L3} that the  conformal structure is $C^1$, it suffices to assume that $f$ is  
$C^2$.
\vskip.1cm

The question whether conformality of $g$ and conjugacy of the periodic data of $f$ and $g$ implies conformality of $f$
is difficult. This is always true for two-dimensional $E^s$ and $E^u $ \cite{KS09}, but  in higher dimensions,
current results have additional assumptions.  In \cite{L2} it is  that for every periodic point $p=g^n(p)$, the restrictions $Dg^n|E^s(p)$ and $Dg^n|E^u(p)$ are scalar multiples of the identity map.
More generally, one can require that the conjugating maps can be chosen uniformly bounded in $p$, 
see Theorem \ref{conf2}  below.

\subsection{Counterexamples and irreducibility} \label{example}

The general higher dimensional case is much more complicated.
In \cite{L1} de la Llave constructed examples showing that conjugacy of the periodic data is not sufficient in general for $C^1$ smoothness of $H$, see also \cite{G}. 
In these examples, $L$ is an automorphism of $T^4$ of the form
$$
L(x,y)=(Ax, By), \quad (x,y)\in \T^2\times \T^2,\
 $$
 where matrices $A,B\in SL(2,\R)$ have eigenvalues $\lambda, \lambda^{-1}$ and 
 $\mu, \mu^{-1}$, respectively, with $\mu> \lambda>1.$ 
 Let $v$ be an eigenvector of $A$ corresponding to $\lambda$, and let $f$ be a perturbation of $L$ of the form 
$$
  f(x,y)= (Ax+ \varphi(y)v,\, By).
$$
The maps $L$ and $f$ have conjugate periodic data since $L^n$ and $D_pf^n$ have the same distinct eigenvalues. 
A conjugacy $H$ is found in the form 
$$\,H(x,y)=(x+\psi(y)v,\, y).\,
$$
This yields a cohomological equation $\varphi(y)+\psi(By)=\lambda \psi(y)$
with a continuous solution
$$
\psi(y)= \lambda^{-1} \sum_{k\ge 0} \lambda^{-k} \varphi(B^ky).
$$
If $\varphi$ and hence $\psi$ are $C^0$ small, then $H$ is the unique conjugacy close to the identity.
However, $\psi$ is usually not differentiable in $y$, for example, for 
$\varphi(y_1,y_2)=\varepsilon \sin (2\pi y_1)$.
\vskip.1cm 
Extending these examples, Gogolev showed in \cite{G} that, outside of the conformal case, 
 irreducibility of $L$ is a necessary 
assumption for conjugate periodic data to imply smoothness of $H$. 
We recall that $L$ is {\em irreducible}\, if it has no nontrivial rational invariant subspace or, equivalently, if its characteristic polynomial is irreducible over $\Q$. While the eigenvalues of an irreducible $L$ are always simple, different  eigenvalues may have the same absolute value, and hence the Lyapunov exponents of $L$ are not necessarily simple.

\vskip.2cm

\subsection{Irreducible $L$ with simple Lyapunov spectrum.}
We now focus on irreducible $L$ and discuss results and approaches for this case. 
We consider the splittings into Lyapunov subspaces  for $L$ and $f$.
Let $1<\rho_1 <\dots <\rho_\ell$ be the distinct moduli of the unstable eigenvalues of $L$ and let 
$$E^{u,L}= E_1^L \oplus E_2^L \oplus \dots \oplus E_\ell^L$$ 
be the corresponding splitting of $E^{u,L}.$
Since $f$ is $C^1$-close to $L$, its unstable sub-bundle splits into a direct sum of $\ell$ invariant 
H\"older continuous sub-bundles close to the ones for $L$ \cite[Section 3.3]{Pes}:
 $$  
 E^{u,f}= E_1^f \oplus E_2^f \oplus \dots \oplus E_\ell^f.
 $$ 
They integrate to H\"older foliations $W_{i}^f$ with $C^{1+\text{H\"older}}$ leaves. While this can be seen using partially hyperbolic theory, unique integrability was obtained directly in \cite{G}. So one can try showing that
$H$ is $C^{1+\text{H\"older}}$ along these foliations. Once this is done, the regularity of $H$ on $\T^d$
follows from Journ\'e lemma.

Local rigidity results for irreducible $L$ with simple Lyapunov exponents were obtained by Gogolev and Guysinski in \cite{GG,G}. 
Even though  in this case  each $E_i^L$ is one-dimensional, there are major difficulties in 
extending Theorem \ref{mainL1}. 
One is showing that $H(W_{i}^f)=W_{i}^L$ for each $i$. 
This always holds for the full unstable foliation and for $i=1$, but it is not true for other foliations in general. Indeed, while the weak flag for $f$ is always mapped by $H$ to that of $L$ \cite{G}, it may not be the case
 for the strong flags. In \cite{G}, the property $H(W_{i}^f)=W_{i}^L$ was proved by an intricate inductive process, 
 which uses partial smoothness of $H$ obtained in the previous steps. The argument also relies on density 
 of leaves of these foliations, which follows from irreducibility of $L$. 

These arguments are difficult to extend to an Anosov diffeomorphism $g$ in place of $L$, as there 
may be no continuous invariant splitting of $W^{u,g}$.
However, it was proved that conjugacy of the periodic data implies that $H$ is $C^{1+\text{H\"older}}$
in the case when $g$ is sufficiently close to an irreducible $L$ with simple spectrum: 
for $\T^3$ in \cite{GG} and for higher dimensions in  \cite{G} under an additional assumptions 
on density of leaves.  Another difficulty in this case is that, unlike $W^{u,g}$, the foliations $W_{i}^g$ are not 
absolutely continuous in general. Apart from these results, as well as the two-dimensional and 
conformal ones above, the problem of smoothness of conjugacy between an arbitrary Anosov 
 diffeomorphism $g$ and its perturbation is open. The results in \cite{G,SY}  for irreducible $L$ 
 with simple Lyapunov exponents were extended by DeWitt in \cite{dW} to the case when $L$ is an Anosov
  automorphisms of a nilmanifold. This is a more complicated setting due to interactions between
  dynamics and the nilpotent structure, so the notion of irreducibility has to be appropriately modified 
  and an extra assumption of sorted spectrum has to be added.

\vskip.1cm
The following theorem generalized the result for toral automorphisms in \cite{G}. 
It allows pairs of eigenvalues of the same modulus. Unlike automorphisms with simple spectrum,  
toral automorphisms $L$ in Theorems \ref{main11}  are generic in the sense that
the proportion of matrices of norm at most $T$  in $SL(d,\Z)$ satisfying the assumptions 
tends to 1 as $T \to \infty$.

\begin{theorem} \cite{GKS11}\label{main11}
Let $L$ be an irreducible Anosov automorphism  of $\T^d$ such that no three of its eigenvalues 
have the same modulus. Let  $f$ be a $C^2$ diffeomorphism of $\T^d$ sufficiently $C^1$-close to $L$. 
If the derivative $D_pf^n$ is conjugate to $L^n$ whenever $p=f^np$, 
then $f$ is $C^{1+\text{H\"older}}$ conjugate to $L$. 
\end{theorem}

\noindent Using similar arguments, outlined below, we obtain another version of this theorem.

\begin{theorem} \label{main11b} 
Let $L$ be an irreducible Anosov automorphism of $\T^d$ and let $f$ be a $C^2$ diffeomorphism of $\T^d$ sufficiently $C^1$-close to $L$.
Suppose that for each $p=f^np$ there is a matrix $C_p$ so that $D_pf^n\,= \,  C_p^{-1} \circ L^n \circ C_p$, 
and the set of all matrices $C_p$ is bounded in $GL(d,\R)$.
Then $f$ is $C^{1+\text{H\"older}}$ conjugate to $L$.

\end{theorem}

\subsection{Rigidity of the Lyapunov spectrum}

The results in \cite{G} for irreducible $L$ with simple Lyapunov exponents were extended
by Saghin and Young in \cite{SY} in a different direction: the assumption of conjugacy of the 
periodic data was weakened to equality of Lyapunov exponents of $L$ to those of $f$ with 
respect to the $f$-invariant volume $\mu$. They used the notion of leaf-wise entropy and 
a leaf-wise analog of Ledrappier's results in \cite {L} to show that $\mu$ has absolutely continuous
conditional measures on the leaves of foliation $W_{i}^f$, which are mapped by $H$ to the 
corresponding conditional measures for $L$. Then H\"older continuity of the densities of these  conditional measures implies, as in the arguments for Theorem \ref{mainL1}, that $H$ is is $C^{1+\text{H\"older}}$ 
along $W_{i}^f$. The following result established in \cite{GKS20} generalized this to  a 
broader class of toral automorphisms, which is generic in the same sense as for Theorem \ref{main11}.


\begin{theorem} \cite{GKS20}\label{main20}
Let $L$ be an irreducible  Anosov automorphism of $\T^d$  such that no three of 
its eigenvalues have the same modulus and $L$ has no pairs of eigenvalues
of the form $\lambda, -\lambda$ or $i\lambda, -i\lambda$, where $\lambda$ is real.
Let  $f$ be a volume-preserving $C^2$ diffeomorphism of $\T^d$ sufficiently 
$C^1$-close to $L$. If the Lyapunov exponents of $f$ with respect to the volume 
are the same as the Lyapunov exponents of $L$, then $f$ is $C^{1+\text{H\"older}}$ conjugate to $L$.\end{theorem}

It may be surprising that equality of exponents for a single measure yields smoothness of the conjugacy, 
but this result is specific to linear $L$. It relies on the property that the measure of maximal entropy for 
$L$ coincides with the volume. The equality of the Lyapunov exponents implies that the same property holds for $f$. Below we deduce a different version of the theorem. The argument showcases the source of rigidity. 

\begin{corollary} \label{main20mme}
Let $L$ be as in Theorem \ref{main20} and let  $f$ be a $C^2$ diffeomorphism of 
$\T^d$ sufficiently $C^1$-close to $L$. 
If the Lyapunov exponents of $f$ with respect to its measure 
of maximal entropy are the same are the same as the Lyapunov exponents of $L$, 
then $f$ is $C^{1+\text{H\"older}}$ conjugate to $L$.

\end{corollary}
We denote by $m$ the Lebesgue measure on $\T^d$. Then $\mu:=(H^{-1})_*(m)$ is the measure of maximal entropy for $f$
 since 
 $$
\mathbf{h}_{top}(f)= \mathbf{h}_{top}(L) = \mathbf{h}_{m}(L)  =\mathbf{h}_{\mu}(f).
$$
We denote the Lyapunov exponents of $f$ with respect to $\mu$ by $\lambda^f$,
and the Lyapunov exponents of $L$ by $\lambda^L$. If these exponents  are the same 
 then by Ruelle inequality we have
 $$
\mathbf{h}_{top}(f)=\mathbf{h}_{\mu}(f)\le  \sum_{\lambda ^f >0} \lambda ^f =  \sum_{\lambda ^L >0} \lambda ^L = \mathbf{h}_{m}(L) = \mathbf{h}_{top}(L) =\mathbf{h}_{top}(f).
$$
Thus equality holds in Ruelle inequality, which implies by \cite{L} that the conditional measures 
of $\mu$ on $W^{u,f}$ are absolutely continuous. The same argument with $f^{-1}$ shows that
 the conditional measures of $\mu$ on $W^{s,f}$ are also absolutely continuous. 
 By the local product structure of the measure of maximal entropy, $\mu$ is absolutely continuous.
 Hence $\mu$ is an $f$-invariant volume 
 and Theorem \ref{main20} applies.


\subsection{Higher-dimensional foliations and linear cocycles}

Now we  describe the techniques and results used in Theorems \ref{main11}, \ref{main11b}, 
and \ref{main20} to deal with higher-dimensional sub-bundles $E_i^L$ and $E_i^f$. 
Suppose that  dim$\,E_i^L=$ dim$\,E_i^f$ is at least two.  Even if  $H$ maps $W_{i}^f$ to $W_{i}^L$, and absolutely continuous conditional measures on the leaves of foliation $W_{i}^f$ to the 
corresponding conditional measures for $L$, smoothness of $f$ along $W_i^f$ does not follow. 
The key new step in Theorems \ref{main11}, \ref{main11b}, and \ref{main20} was to establish 
conformality of the derivative cocycle $Df|E_i^f$. It used results on linear cocycles
over hyperbolic systems, which we now discuss.

\vskip.1cm 

Let $f$ be an Anosov diffeomorphism of $\M$ and let
 $A$ be a map from $\M$ to  $GL(m,\R)$.
The $GL(m,\R)${\em -valued  cocycle over $f$ generated by }$A$
is the map $\A:\,\M \times \Z \,\to GL(m,\R)$ given  by $\,\A(x,0)=\Id\,$ and for $n\in \N$,
$$
 \A(x,n)=\A_x^n = A(f^{n-1} x) \cdots  A(x) \;\,\text{ and }\;\,
\A(x,-n)=\A_x^{-n}= (\A_{f^{-n} x}^n)^{-1}.
$$
The regularity of  $\A$ is defined as that of  if its generator $A$.

A prime example of a linear cocycle is the derivative cocycle.
For $T\M= \M\times \R^m$, one can take $ A(x)=D_xf  \in GL(m,\R),$  and then
$ \A_x^n=D_xf^n.$\,   Similarly, one can consider $A(x)=Df|\E(x)$, where $\E$ is a $Df$-invariant sub-bundle
 such as $E^s$, $E^u$, or $E_i^f$.
 
We say that cocycles $\A$ and $\B$ are (measurably or continuously) hohomologous if there exists a (measurable or continuous) function $\c:\M\to GL(m,\R)$ such that
$$
  \A_x=\c(fx)\,\B_x\,\c(x)^{-1} \quad\text{for all }x\in \M.
$$
The function $C$ is called a {\em conjugacy} or a {\em transfer map} between $\A$ and $\B$.
\vskip.1cm
The next theorem established conformality of a cocycle based on its periodic data.

\begin{theorem} \cite{KS10} \label{conf2}
Let $\A$ be a $GL(m,\R)$-valued H\"older continuous cocycle over
a transitive Anosov diffeomorphism $f:\M\to\M$.
Suppose that for each periodic point $p=f^n(p)$ in $\M$ there exists a matrix $C_p$ such that 
$\,C_p^{-1}\,\A_p^n\, C_p\,$ is conformal.

If either $m=2$, or $m>2$ and the set of matrices $C_p$  is bounded in $GL(m,\R)$,
then $\A$ is H\"older continuously cohomologous to a {\em conformal cocycle},
i.e., a cocycle with values in the conformal subgroup.
\end{theorem}

For $m>2$, the result does not hold without the boundedness assumption \cite{KS10}.
This theorem yields conformality of the derivative $Df$ on each two-dimensional $E_i^f$
in Theorems \ref{main11} and \ref{main11b}, and explains the difference in their assumptions.
Theorem \ref{main20} uses the continuous amenable reduction results \cite{KS13}, 
which for $m=2$ yield that a cocycle is cohomologous to conformal  if it does not have a 
continuous invariant sub-bundle or a continuous invariant field of two lines \cite{GKS20}. 
 If $L$ has no pairs of eigenvalues $\lambda, -\lambda$ or $i\lambda, -i\lambda$ with $\la \in \R$, 
 then $L|E_i^L$ has no such invariant objects, and $Df | E_i^f$ for a small perturbation has the same property.

\vskip.1cm
Once it is established that $H(W_{i}^f)=W_{i}^L$, regularity of $H$ along $W_{i}^f$ is proved using
conformality of $L$ on $E_i^L$ and of $Df$ on $E_i^f$. After conjugating $L|E_i^L$ and $Df | E_i^f(x)$ 
to conformal cocycles, the norms give scalar cocycles $a(x)=\rho_i =\|L|E_i^L\|$  and $b(x)=\|Df | E_i^f(x)\|$.  
In Theorems \ref{main11} and \ref{main11b}, these H\"older continuous scalar cocycles have equal periodic data, and hence they are H\"older continuously cohomologous by the Liv\v{s}ic periodic point theorem  \cite{Liv2}. In Theorem \ref{main20}, there is a measurable conjugacy between $a$  and $b$, 
 which is obtained from the Jacobian of $H$ along $W_{i}^f$, and it follows by the measurable Liv\v{s}ic theorem  \cite{Liv2} that the conjugacy is H\"older continuous.
 
Continuous conjugacy  between $a$  and $b$ implies that the ratio of norms $\|L^n|E_i^L\|$ and $\|Df^n | E_i^f\|$ 
is bounded above and below uniformly in $x$ and $n$. This is used to show that $H$ is bi-Lipschitz continuous along $W_{i}^f$, which yields differentiability of $H$ almost everywhere on each leaf. 
The derivative $DH$ is then a bounded measurable conjugacy between linear cocycles
$Df|E_i^f$ and $L|E_i^L$. Thus to conclude that $H$ is $C^{1+\text{H\"older}}$ along $W_{i}^f$
it suffices to establish H\"older continuity of the conjugacy $DH$. 
\vskip.1cm

Continuity of a measurable conjugacy between H\"older continuous cocycles over hyperbolic systems has been extensively studied. 
It always holds for scalar cocycles \cite{Liv2}, but may fail
already for $GL(2,\R)$-valued cocycles with more than one Lyapunov exponent, even when both 
generators are close to the identity \cite{PW01}. Continuity was obtained under various compactness, boundedness, and conformality assumptions on both or one of the cocycles \cite{PaP97,Pa99,Sch99,S13,S15,B}. Results for conformal cocycles were used to obtain H\"older continuity of $DH|E_i^f$ above.


\subsection{New results for cocycles and smoothness of a weakly differentiable conjugacy.}
In our recent work \cite{KSW} we establish H\"older continuity of a measurable conjugacy 
in the optimal setting of  cocycles with one Lyapunov exponent.

\begin{theorem}  \cite{KSW}\label{measurable conjugacy}
Let $f:\M\to \M$ be a  transitive $C^{1+\text{H\"older}}$   Anosov diffeomorphism,
and let $\A$ and $\B$ be $\beta$-H\"older linear cocycles over $f$.
Let  $\mu$ be an ergodic $f$-invariant  measure on $\M$ with full support and local product structure.

Suppose that $\A$ has one Lyapunov exponent at every periodic point and $\B$ is fiber bunched.
Then any $\mu$-measurable conjugacy between $\A$ and $\B$ is $\beta$-H\"older continuous, i.e.,
coincides with a $\beta$-H\"older continuous conjugacy on a set of full measure.

\end{theorem}

Fiber bunching is a technical condition weaker than having one exponent.
It means that non-conformality of the cocycle is dominated by the contraction and expansion in the base.
In the proof we apply continuous amenable reduction \cite{KS13} to  $\A$ and obtain a block-triangular structure with conformal blocks on the diagonal. Using the measurable conjugacy $\c$, we obtain a similar measurable block structure for $\B$ and establish its continuity. Then continuity of  the diagonal blocks of $\c$ follows from \cite{S15}. To prove continuity of the off-diagonal blocks of $\c$ we use an inductive process. For this we establish a result on continuity of a measurable conjugacy for vector-valued 
cocycles twisted by a bounded linear cocycle. 

As a corollary, we obtain the following result for perturbations of constant cocycles,
which we then use in the new local rigidity results.

\begin{theorem} \cite{KSW}\label{constant cocycle}
Let $f$ and $\mu$ be as in Theorem \ref{measurable conjugacy} and
let  $\A$ is be a constant $GL(m,\R)$-valued cocycle over $f$. Then
for any H\"older continuous $GL(m,\R)$-valued cocycle $\B$ sufficiently $C^0$ close to $\A$, any $\mu$-measurable conjugacy between $\A$ and $\B$ is  H\"older continuous.
\end{theorem}
\noindent Also, we obtain estimates of the H\"older exponent and H\"older constant of the conjugacy.

\vskip.2cm
As we observed, existence of a continuous conjugacy between the derivative cocycles is closely related to 
smoothness of $H$. However, the relationship is not straightforward. If $Df$ is continuously conjugate to 
$L$ and $L$ is irreducible, then Theorem \ref{main11b} yields that $H$ is $C^{1+\text{H\"older}}$. 
Without irreducibility,  however, existence of {\it some} continuous conjugacy between the derivative 
cocycles $Df$ and 
$L$ does not imply in general that $H$ is $C^1$. In fact, if all eigenvalues of $L$ are simple with distinct moduli, then conjugacy of $D_pf^n$ and $L^n$ whenever $p=f^n(p)$ always gives H\"older conjugacy of the cocycles
(as the cohomological equations splits into scalar ones for the restrictions to $E_i$),  but $H$ may not be $C^1$ as the counterexample above shows.

\vskip.2cm

However, if $H$ is differentiable in a weak sense, we show that $H$ is $C^{1+\text{H\"older}}$.
This result holds for an
arbitrary hyperbolic automorphism without any irreducibility assumption.
We denote by $W^{1,q}(\T^d)$ the Sobolev space of $L^q$
functions with $L^q$ weak partial derivatives of first order. Note that any Lipschitz function is in $W^{1,\infty}(\T^d)$. 


While $H$ satisfying \eqref{IntrConj} is not unique,  there is a unique {\em conjugacy $C^0$ close to the identity}. This is the unique $H$ in the homotopy class of the identity with $H(p)=0$, where $p$ is the fixed point of $f$ closest to $0$.

\begin{theorem} \label{HolderConjugacy}
Let $L$ be a hyperbolic automorphism of $\T^d$ and let $f$ be a $C^{1+\text{H\"older}}$ diffeomorphism of $\T^d$ which is $C^1$ close to $L$. Suppose that for some conjugacy $H$
between $f$ and $L$, either $H$ or $H^{-1}$ is in  $W^{1,q}(\T^d)$ with $q>d$.
Then $H$ is a $C^{1+\text{H\"older}}$ diffeomorphism.
\vskip.1cm

More precisely, there is a constant $\beta_0=\beta_0(L)$, $0<\beta_0\le1$, so that for any
 $0<\beta' <\beta_0$ there exist constants $\delta>0$ and $K>0$ such that for any $0<\beta \le \beta'$
the following holds.

For any ${C^{1+\beta}}$ diffeomorphism $f$ with $\|{f-L}\|_{C^{1}}<\delta$, if some conjugacy
between $L$ and $f$, or its inverse, is in $W^{1,q}(\T^d)$, $q>d$, then any conjugacy  is  a
$C^{1+\beta}$ diffeomorphism. Moreover, for the conjugacy $H$ that is $C^0$ close to the identity,
\begin{equation} \label{C1H est}
\|{H-I}\|_{C^{1+\beta}}\leq K\|{f-L}\|_{C^{1+\beta}}.
\end{equation}
\end{theorem}

In the proof we differentiate the conjugacy equation  $L\circ H= H \circ f$ to obtain 
\begin{equation} \label{IntrConjD}
L\circ DH= DH \circ Df,
\end{equation}
where $DH$ is the Jacoby matrix of weak partial derivatives. We use the assumption that
 $H$ is in $W^{1,q}(\T^d)$ with $q>d$ to show that $f$ preserves an absolutely continuous measure $\mu$,
  that $DH$ gives the differential of $H$ $\mu$-a.e on $\T^d$, that  \eqref{IntrConjD} holds $\mu$-a.e, 
  and that $DH$ is invertible $\mu$-a.e. Hence \eqref{IntrConjD} shows that $DH$ is a measurable
   conjugacy between $L$ and $Df$. 
Then Theorem \ref{constant cocycle} yields that $DH$ is H\"older continuous, and so $H$ is a $C^{1+\text{H\"older}}$ diffeomorphism.

The inequality \eqref{C1H est} is obtained using the estimate for the conjugacy between cocycles in Theorem \ref{constant cocycle}. It plays an important role in establishing higher regularity of $H$ in Theorem \ref{th:4} below.

\section{Higher regularity  of the conjugacy} 

In dimension two and in the higher-dimensional conformal case, 
Theorems \ref{mainL1} and \ref{conf} yield that if $C^\infty$ Anosov diffeomorphisms $f$ and $g$ 
 are conjugate by a $C^1$ diffeomorphism $H$, then $H$ is $C^\infty$. 
The general higher dimensional case is much more complicated. The problem of the exact 
regularity of $H$ is subtle: for any $k \in \N$ and any $d\ge 4$ there exists a reducible hyperbolic automorphism $L$ of $\T^d$ and its analytic perturbation $f$ such that the conjugacy $H$ is $C^k$ but is not $C^{k+1}$ \cite{L1}. This was demonstrated by examples as in Section \ref{example} with an 
appropriate relation between the eigenvalues $\lambda $ and $\mu$.
On the other hand, for a given $L$ (or a nonlinear $g$) there is $k(L)$ such that if $H$ is  $C^{k(L)}$ then it is $C^{\infty}$ \cite{L1}. 

Theorems \ref{main11}, \ref{main11b}, and \ref{main20} yield only that $H$ is $C^{1+\text{H\"older}}$.
 This low smoothness is due to the method of the proof, which shows regularity
of $H$ along foliations $W^i_f$, whose leaves are typically only $C^{1+\text{H\"older}}$ smooth. 
Nevertheless, Gogolev conjectured \cite{G} that the conjugacy in Theorem \ref{mainL1}
should be $C^{r-\e}$, if $f$ is $C^r$. The only progress in this direction until now was the result of Gogolev 
for automorphisms of $\T^3$ with real spectrum in \cite{G17}.  

\vskip.1cm
In the next theorem we obtain $C^\infty$ smoothness of a conjugacy to a $C^\infty$ perturbation 
$f$ assuming that $L$ is weakly irreducible, which is defined as follows. Let $\R^d=\oplus_{\rho_i} E^i$ be the splitting where $E^i$ is the sum of generalized eigenspaces of $L$ corresponding to the eigenvalues of modulus $\rho_i$, and let $\hat E^i=\oplus_{\rho_j\neq \rho_i} E^j.$
 We say that $L$ is \emph{weakly irreducible} if each $\hat E^i$ contains no nonzero elements of $\Z^d$.
 Irreducibility over $\Q$ implies weak irreducibility. Indeed, if there is a nonzero integer point $n \in \hat E^i$ then $span \{ L^m n : m\in \Z \} \subset \hat E^i$ is a nontrivial rational invariant subspace.
Weak irreducibility is determined by the characteristic polynomial of $L$ as follows.

 \begin{lemma} \cite{KSW} \label{weak irred}
A matrix $L\in GL(d,\Z)$ is weakly irreducible if and only if there is a set $\Delta \subset \R$ so that for each
irreducible over $\Q$ factor of the characteristic polynomial of $L$ the set of moduli of its roots
equals $\Delta$.
\end{lemma}

It follows that if $L$ is irreducible or weakly irreducible then the following
matrices are weakly irreducible
$$
\left(\begin{array}{cc}L & 0 \\ 0 & L \end{array}\right) \quad \text{and}\quad
\left(\begin{array}{cc}L & I \\ 0 & L \end{array}\right).
$$
These matrices are not irreducible and the second one  is not even diagonalizable. 
So while an irreducible $L$ is always conformal (in some metric) on each Lyapunov subspace, 
weakly irreducible $L$ may have Jordan blocks.

\begin{theorem} \cite{KSW}\label{th:4} Let $L$ be a weakly irreducible hyperbolic automorphism of $\T^d$.
Then there is $r=r(L)\in\N$ so that for any $C^\infty$ diffeomorphism $f$   which is $C^r$ close to $L$ the following holds. If for some conjugacy $H$ between $f$ and $L$ either $H$ or $H^{-1}$ is in the Sobolev space $W^{1,q}(\T^d)$ with $q>d$, then any conjugacy between $f$ and $L$ is a $C^\infty$ diffeomorphism. 
\end{theorem}

Our approach is completely different from the previous local rigidity results.  
Theorem \ref{HolderConjugacy} is used as the first step in the proof of Theorem~\ref{th:4}
to obtain $C^{1+\b}$ regularity of $H$ and the estimate \eqref{C1H est}. 
Then to establish $C^{\infty}$ smoothness of $H$ we
use an iterative method which is somewhat similar to the traditional KAM scheme. 
However, KAM is primarily used for elliptic systems and not for hyperbolic ones. 
The main ingredient in the iterative step is obtaining an  approximate $C^\infty$  solution for the 
linearized conjugacy equation. 
\vskip.1cm 
Assuming $f(0)=H(0)=0$, we lift $f$ and $H$ to $\R^d$ as $\tilde H=\Id+h$ and $\tilde f=L+R$, 
where $h, R:\R^d \to \R^d$ are $\Z^d$-periodic functions.
The lifts satisfy the conjugacy equation $L\circ \tilde H=\tilde H\circ \tilde f$, which yields
$L\circ h - h\circ \tilde f=R.$ The latter projects to the torus as
$$
 L\circ h - h\circ  f=R.
$$
This is a twisted by $L$ cohomological equation over $f$ for $\R^d$-valued functions $h$ and $R$ on $\T^d$, where $R$ is $C^{\infty}$ and $h$ is $C^{1+\b}$. The iterative scheme relies on finding an {\em approximate $C^\infty$ solution} with good estimates for the {\em linearized equation} over $L$:
$$
 L\circ h' - h'\circ  L=Q, \quad\text{where}\quad Q-R= h\circ  f-h\circ  L. 
$$
Conjugating $f$ by the $C^\infty$ diffeomorphism $H'=\Id-h'$ we  get a new $C^\infty$ diffeomorphism  $f'$, 
which is much closer to $L$ and is still $C^{1+\b}$ conjugate to $L$. Continuing this iterative process 
and establishing convergence in a suitable smooth topology proves the theorem. 

Closest to our setting, KAM techniques were used in \cite{Damjanovic4} to prove $C^\infty$ local rigidity 
for some $\Z^2$ actions by partially hyperbolic toral automorphisms.  The structure of a higher 
rank action was used there in an essential way to construct such an approximate solution.
In our case these methods do not apply, and instead the argument relies on existence of a
 $C^{1+\b}$ conjugacy with the estimate \eqref{C1H est}.  The linearized equation is analyzed 
 using Fourier coefficients. However, expansion/contraction of the twist $L$ requires working in
 higher smoothness classes, while our functions $h$, and hence $Q$, are only $C^{1+\b}$.
To employ our low regularity data we use the Lyapunov splitting $\R^d= \oplus E^i$  for $L$ and consider projections
$$ L_i \circ h_i-h_i\circ L =Q_i, \quad \text{where } L_i=L|{E^i}.$$
Differentiating along $E^i$, we ``balance" the twist by the derivative
$$
L_i \cdot D_i h_i-(D_i h_i)\circ L \cdot L_i =D_i Q_i \quad
$$ 
which allows us to use the existence of H\"older solution $Dh$.

\vskip.3cm
Applying Theorem \ref{th:4} we improve the regularity of the conjugacy from $C^{1+\text{H\"older}}$ to $C^\infty$ in Theorems \ref{main11} and \ref{main20}.

\begin{corollary} \label{local PD}
Let $L:\T^d\to\T^d$ be an irreducible Anosov automorphism
such that no three of its eigenvalues have the same modulus.
Let  $f$ be a $C^\infty$ diffeomorphism which is $C^r$  close to $L$ such that the
derivative $D_pf^n$ is conjugate to $L^n$ whenever $p=f^n(p)$.
Then $f$ is $C^{\infty}$ conjugate to $L$.
\end{corollary}

\begin{corollary} \label{local LS}
Let $L:\T^d\to\T^d$ be an irreducible Anosov automorphism such that no three of
its eigenvalues have the same modulus and there are no pairs of eigenvalues
of the form $\lambda, -\lambda$ or $i\lambda, -i\lambda$, where $\lambda\in \R$.
Let  $f$ be a volume-preserving $C^\infty$ diffeomorphism of $\T^d$ sufficiently
$C^r$-close to $L$. If the Lyapunov exponents of $f$ with respect to the volume
are the same as the Lyapunov exponents of $L$, then $f$ is $C^{\infty}$ conjugate to $L$.
\end{corollary}


\end{document}